\documentclass[12pt]{article}
\usepackage{amstext}
\usepackage{fancyhdr}
\usepackage{etoolbox}
\usepackage{mathtools}
\usepackage{amstext}
\usepackage{fancyhdr}
\usepackage{amsfonts,graphicx,bezier, amssymb}
\usepackage{amsmath}
\usepackage{caption}
\usepackage{mathtools}
\usepackage[utf8]{inputenc}
\usepackage{tikz}
\usetikzlibrary{arrows.meta}
\newtheorem{thm}{Theorem}[section]
\newtheorem{defn}[thm]{Definition}
\newtheorem{prop}[thm]{Proposition}
\newtheorem{lem}[thm]{Lemma}
\newtheorem{cor}[thm]{Corollary}
 \newtheorem{exam}[thm]{Example}
 \newenvironment{prf}{\noindent{\bf{\it{Proof}:}}~~}{\hfill\rule{1ex}{1ex}\vskip1.5ex}

\newcommand{\Z}{\mathbb Z}
\newcommand{\Q}{\mathbb Q}

\newcommand{\beqa}{\begin{eqnarray}}
\newcommand{\enqa}{\end{eqnarray}}
\newcommand{\beq}{\begin{eqnarray*}}
	\newcommand{\enq}{\end{eqnarray*}}

 \newtheorem{rem}{Remark}[section]

\begin{document}
	\begin{center}
		{\bf\Large Nil modules and the envelope of a submodule}
		
		\vspace*{0.2cm}

		\begin{center}
			
			  	David Ssevviiri\footnote{Department of Mathematics, Makerere University, 
				P.O. BOX 7062, Kampala, Uganda}$^{,}$\footnote{Corresponding author} and Annet Kyomuhangi\footnote{Department of Mathematics, Busitema University,
				P.O. BOX 236, Tororo, Uganda} \\
			E-mail:  david.ssevviiri@mak.ac.ug and annet.kyomuhangi@gmail.com
			
		\end{center}

	\end{center}	
	\vspace*{0.2cm}
\begin{abstract}
		Let $R$ be a commutative unital ring and $N$ be a submodule of an $R$-module $M$.    We  show that: 1) the semiprime radical is an invariant  on submodules  generated by the ascending chain of envelopes of a given submodule; 2) for rings that satisfy the radical formula, $\langle E_M(0)\rangle$ is an idempotent radical leading to a torsion theory whose torsion class has  nil $R$-modules and the torsion-free class  has reduced $R$-modules; and, 3) Noetherian uniserial modules satisfy the semiprime radical formula and their semiprime radical is a nil module.
		
	\end{abstract}
	
	
	\vspace*{0.3cm}
	
	{\bf MSC 2010} Mathematics Subject Classification:   16S90, 13C13 
	\section{Introduction}
	
	\paragraph\noindent
	Throughout the paper, all rings $R$ are commutative with unity and all modules $M$ are left $R$-modules. The notions of a module or a ring that satisfy the radical formula and modules that satisfy the semiprime radical formula are a common theme in the literature, see for instance, \cite{Azizi, Jenkins, Lee2018, Le, McCasland, Parkash, Sharif, PSmith, Steven}. A proper submodule $N$ of an $R$-module $M$ is \emph{prime} if for all $a\in R ~\text{and}~ m\in M$, $am\in N$ implies that either $m\in N$ or $aM\subseteq N.$ A proper submodule $N$ of an $R$-module $M$ is {\it semiprime} if for all $a\in R$ and $m\in M$, $a^2m\in N$ implies that $am\in N$. 	A module is prime (resp. reduced) if its zero submodule is prime (resp. semiprime).  The intersection of all prime (resp. semiprime) submodules of $M$ containing $N$ is called the {\it prime radical} ({\it semiprime radical}) of $N$ which
	 is  denoted by $\beta(N)$  (resp. $S(N)$). If $N=0$, then we write $\beta(M)$ (resp. $S(M)$) and call this the prime radical (resp. semiprime radical) of $M$.
	
	\paragraph\noindent 
	
	For a submodule $N$ of an $R$-module $M$,  the envelope of  $N$ in $M$ is the set
	$$E_M(N):= \left \{am~|~ a^{k}m\in N, ~a\in R,~m\in M,~ \text{for some }k\in \mathbb{Z}^{+} \right\}.$$ 	
	$E_{M}(0)$ was considered as the module analogue of $\mathcal{N}(R),$ the collection of all nilpotent elements of the ring $R,$ since for  any ring $R$, $E_{R}(0)= \mathcal{N}(R).$ In general, $E_M(N)$ is not a submodule of $M$. The submodule $\langle E_M(N)\rangle$ generated by $E_M(N)$ is the module analogue of the radical of an ideal $I$ of $R$ given by $\sqrt{I}:=\{a~|~a^{k}\in I \text{ for some }k\in \Z^{+}\}.$
	A module $M$ {\it satisfies the radical formula}  (resp. {\it satisfies the semiprime radical formula})  if for every submodule $N$ of $M$, $\langle E_M(N)\rangle=\beta(N)$ (resp. $\langle E_M(N)\rangle=S(N)$). A ring $R$ {\it satisfies the radical formula} if every $R$-module satisfies the radical formula.
	
	\paragraph\noindent
	
	In Section 2, we define what we call the envelope functor. It associates to an $R$-module $M$, a submodule $\langle E_M(0)\rangle$. We show that this functor is a preradical [Proposition \ref{fun}] which is not a radical  in general [Example \ref{ex}].	 	As in \cite{Azizi} and \cite{Lee2018}, if $N$ is a submodule of $M$, we  define $E_0(N):=N$, $E_1(N):=E_M(N)$,  
	$E_2(N):= E_M(\langle E_{1}(N)\rangle)$, $\cdots$, $E_n(N):= E_M(\langle E_{n-1}(N)\rangle)$ which forms an ascending chain $E_0(N)\subseteq E_1(N)\subseteq E_3(N)\subseteq \cdots$ of envelopes. The first main theorem, Theorem \ref{chain} says that if the chain of envelopes above terminates at $E_n(N)$, then for all $i=0,1,2,\cdots, n$; $S(\langle E_{i}(N)\rangle)= \langle E_{n}(N)\rangle$, i.e., the semiprime radical is an invariant on the submodules generated by the envelopes in the finite chain.
	
	\paragraph\noindent
	In Section 3, we introduce nil modules and show that $\langle E_{M}(0)\rangle$ is the largest nil submodule of $M$, [Proposition \ref{e}]. The second main theorem, Theorem \ref{so} says that; for rings that satisfy the radical formula, the envelope functor is an idempotent radical and it induces a torsion theory with all the nil $R$-modules as the torsion class and all the reduced $R$-modules as the torsionfree class. We further show that a Noetherian uniserial $R$-module satisfies the semiprime radical formula [Theorem \ref{rad}] and its semiprime radical $S(M)$ is a nil module [Proposition \ref{nil}].

	\section{The envelope functor}
 
	\paragraph\noindent
		A functor $F:R\text{-Mod} \to R\text{-Mod} $ is a {\emph{preradical}} if
		for every $R$-module $M$, $F(M)$ is a submodule of $M$ and for every $R$-homomorphism $\gamma:M \to N,~ \gamma(F(M))\subseteq F(N)$.
		$F $ is a {\emph{radical}} if it is a preradical and for all $M \in R$-Mod,  $F(M/F(M))=0.$

	\begin{prop}\label{fun}
		Let $R$ be a ring and $ M\in R$-Mod.  The endofunctor $$F: R\text{-Mod}\to R\text{-Mod},~M \mapsto \langle E_{M}({0}) \rangle$$ is a preradical.
	\end{prop}
\begin{prf}
	Let $M,N \in R$-Mod and $\gamma: M\to N$ be an $R$-module homomorphism. If $x\in \gamma(F(M))=\gamma(\langle E_{M}(0) \rangle),$ then $x=\gamma(m)$ for some $m\in \langle E_{M}(0) \rangle.$ So, $m=\sum \limits_{i=1}^{t}a_{i}m_{i}$ for some $a_{i}\in R$ and $m_{i}\in E_{M}(0).$ This implies that $m_{i}=b_{i}s_{i}$ and $b_{i}^{k_{i}}s_{i}=0$ for some $s_{i}\in M$, $b_{i} \in R$ and $k_{i}\in \Z^{+}.$ As such, $x= \sum \limits_{i=1}^{t}a_{i}b_{i}\gamma(s_{i}).$ Since $b_{i}^{k_{i}}s_{i}=0,$  $b_{i}^{k_{i}}\gamma(s_{i})=0$  and $b_{i}\gamma(s_{i})\in E_{N}(0).$ Thus, $ x\in \langle E_{N}(0) \rangle=F(N)$ and  $\gamma(F(M))\subseteq F(N).$
\end{prf}

\paragraph\noindent

The functor $F$ in Proposition \ref{fun} will be called the {\it envelope functor}.

\paragraph\noindent 
Let $N$ be a submodule of $M$.  Define  $E_{0}(N):= N$, $E_{1}(N):= E_M(N)$,
	$E_{2}(N):= E_M(\langle E_{1}(N)  \rangle)$, 
	$E_{3}(N):= E_M(\langle E_{2}(N) \rangle)$
	$\cdots$,
	$E_{n}(N):= E_M(\langle E_{n-1}(N) \rangle). $
We get an ascending chain

 \begin{equation}\label{chain1}
 E_0(N)\subseteq E_{1}(N)\subseteq E_{2}(N)\subseteq \cdots \subseteq E_{n}(N) \subseteq \cdots
 \end{equation}
 
of envelopes of submodules of $M.$ $E_{n}(N)$ is  called the {\it $n$th envelope} of $N.$  Chain (\ref{chain1}) first appeared in \cite{Azizi} where modules that satisfy the radical formula of degree $n$ were studied. Let $n\in \Z^{+}\cup \{0\}.$ A module {\it satisfies the radical formula of degree $n$} if $ \langle E_{n}(N)  \rangle=$  $\beta(N)$ for every submodule $N$ of $M$.   
Thus, a module satisfies the radical formula  if it satisfies the radical formula of degree one.

\begin{paragraph}\noindent 

Proposition \ref{env} below shows that a quotient module by the consecutive submodules $\langle E_{n-1}(N)\rangle$ and $\langle E_{n}(N)\rangle$  in the  chain (\ref{chain2}) below

\begin{equation}\label{chain2}
 N\subseteq \langle E_1(N)\rangle \subseteq \langle E_2(N)\rangle \subseteq \cdots
 \end{equation}
 
 coincides with the submodule generated by the envelope of the  quotient $R$-module $M/\langle E_{n-1}(N)\rangle$. Compare with \cite[Lemma 2.1]{AppII}.
\end{paragraph}

\begin{prop}\label{env}
	For any submodule $N$ of an $R$-module $M$  such that $E_{1}(N)\subseteq E_{2}(N)\subseteq \cdots \subseteq E_{n}(N) \subseteq \cdots$ is an ascending chain as defined in $(1),$ we have 
	
$$\frac{\langle E_n(N)\rangle}{	\langle E_{n-1}(N)\rangle}= \langle E_{\frac{M}{\langle E_{n-1}(N)\rangle}}(\bar{0})\rangle.$$

\end{prop}

\begin{prf}
Let $\displaystyle{\bar{x}\in \frac{ \langle E_{n}(N)  \rangle}{ \langle E_{n-1}(N) \rangle}}$, $\bar{x}=\sum \limits_{i=1}^{n}a_{i}x_{i} + \langle E_{n-1}(N) \rangle$ where $a
_{i}\in R$ and $x_{i}\in E_{n}(N).$ By definition the  of $E_n(N)$, $x_{i}\in E_M(\langle E_{n-1}(N) \rangle).$ So, $x_{i}=b_{i}m_{i}$ and $b_{i}^{k_i}m_{i}\in \langle E_{n-1}(N) \rangle $ for $m_{i}\in M, ~b_{i}\in R$ and $k_{i}\in \Z^{+}.$ It follows that, $\bar{x}=\sum \limits_{i=1}^{n}a_{i}b_{i}m_{i} + \langle E_{n-1}(N) \rangle$ and $b_{i}^{k_{i}}m_{i}\in \langle E_{n-1}(N) \rangle $ for $m_{i}\in M, ~b_{i}\in R$ and $k_{i}\in \Z^{+}.$ Thus, $\bar{x}\in \langle E_{\frac{M}{\langle E_{n-1}(N) \rangle }}(\bar{0}) \rangle.$  
Conversely, suppose that $\bar{m}\in \langle E_{\frac{M}{\langle E_{n-1}(N) \rangle}}(\bar{0}) \rangle$, $\bar{m}=\sum_{i=1}^na_i(m_i+\langle E_{n-1}(N) \rangle)$, where $a_i\in R$ and $m_i+\langle E_{n-1}(N) \rangle\in E_{\frac{M}{\langle E_{n-1}(N) \rangle}}(\bar{0})$. By the definition of $E_{\frac{M}{\langle E_{n-1}(N) \rangle}}(\bar{0})$, $m_i+\langle E_{n-1}(N) \rangle=r_i(n_i+\langle E_{n-1}(N) \rangle)$ and $r_i^{k_i}(n_i+\langle E_{n-1}(N) \rangle)\in \langle E_{n-1}(N) \rangle$ for some $n_i\in M$, $r_i\in R$ and $k_i\in \Z^+$ for all $i\in \{1, \cdots, n\}$. So, $\bar{m}= \sum_{i=1}^na_i(r_in_i+\langle E_{n-1}(N) \rangle)= (\sum_{i=1}^na_ir_in_i)+\langle E_{n-1}(N) \rangle$. However, $r_i^{k_i}n_i\in \langle E_{n-1}(N) \rangle$ for all $i\in\{1, \cdots, n\}$ implies that  $r_in_i\in E_M(\langle E_{n-1}(N) \rangle) = E_n(N)$ for all $i\in \{1, \cdots, n\}$. This shows that  $\bar{m}\in \frac{\langle E_{n}(N)\rangle}{\langle E_{n-1}(N) \rangle}$.
\end{prf}

\begin{lem}\label{ksem}
	If $N$ is a submodule of $M$ and $K$ is a semiprime submodule of $M$ such that $N \subseteq K,$ then $$N \subseteq E_M(N) \subseteq K.$$
\end{lem}

\begin{prf}
	If $x \in E_M(N)$, then $x=am$ and $a^{k}m\in N$ for some $m\in M$ and  $k\in \Z^{+}.$ So, $a^{k}m \in K$ since $N\subseteq K.$ By the definition of a semiprime submodule, $x=am \in K.$ The inclusion $N \subseteq E_M(N)$ is well known.
\end{prf}

\begin{lem}\label{small}
	Let $N$ be a submodule of $M$. If $\langle E_M(N)  \rangle$ is a semiprime submodule of $M$, then it is the smallest semiprime submodule of $M$ containing $N,$ i.e., $\langle E_M(N)  \rangle =S(N)$ the semiprime radical of $N$. 
\end{lem}

\begin{prf}
	Suppose there exists a semiprime submodule  $K$ of $M$ such that $N \subseteq K \subseteq E_M(N).$ By Lemma \ref{ksem}, $N \subseteq    E_M(N)  \subseteq K$.	So, $\langle E_M(N)  \rangle\subseteq K$ and $K=\langle E_M(N)  \rangle.$
\end{prf}

\ 
 
 \begin{thm}\label{s}{\rm \cite[Theorem 3.8]{Lee2018}} If $N$ is a submodule of an $R$-module $M$, then 
 $$S(N)= \lim_{\stackrel{\longrightarrow}{n}} \langle E_n(N)\rangle.$$
 \end{thm}

Lemma \ref{for} given below  is well known, see for instance \cite[Remark 2.6(4)]{Steven}.
\begin{lem}\label{for}
	Let $N$ be a submodule of $M$. $E_M(N)=N$ if and only if $N$ is a semiprime submodule of $M.$ In particular, $M$ is a  reduced  $R$-module if and only if $E_{M}(0)=0.$
\end{lem}

\begin{cor}\label{cc}If Chain (\ref{chain2}) terminates at $\langle E_n(N)\rangle$, then $$S(N)= \langle E_n(N)\rangle.$$

\end{cor}
 
\begin{thm}\label{chain}
Let $N$ be a submodule of $M$ and consider Chain (\ref{chain2}).   If this chain terminates at the submodule $\langle E_{n}(N) \rangle,$ then for all $i=0, 1,2, \cdots, n$ we  have $$ S(\langle E_{i}(N) \rangle)=\langle E_{n}(N) \rangle. $$
\end{thm}

\begin{prf}
	Define
	 $F_{1}: =  \langle E_{1}(N) \rangle ,
		F_{2}$ as the submodule of $ M$ such that
			$F_{2}/F_{1}:=  \langle E_{M/F_{1}}(\bar{0}) \rangle,   \cdots,$ and  $F_n$ the submodule of $M$ such that 
		$F_{n}/F_{n-1}:=  \langle E_{M/F_{n-1}}(\bar{0}) \rangle $. We get an ascending chain
\begin{equation}
N \subseteq F_{1} \subseteq F_{2} \subseteq \cdots \subseteq F_{n} \subseteq \cdots
\end{equation} of submodules of $M.$ Chain $(2)$ terminates at $F_{n}$ if and only if $ \langle E_{M/F_{n}}(\bar{0}) \rangle =\bar{0}$ if and only if $M/F_{n}$ is a reduced $R$-module (Lemma \ref{for}), i.e., if and only if $F_{n}$ is a semiprime submodule of $M.$ By Lemma \ref{env},   $F_{n}= \langle E_{n}(N)  \rangle.$  
Since $\langle E_n(N)\rangle$ is a semiprime submodule of $M$ by Proposition \ref{small},

\begin{equation}\label{eq2}
\langle E_n(N)\rangle =\langle E_M(\langle E_{n-1}(N)\rangle )\rangle =  S(\langle E_{n-1}(N)\rangle).
\end{equation}

By Corollary \ref{cc}, 
\begin{equation}\label{eq3}
\langle E_n(N)\rangle= S(N).
\end{equation}

From equalities   (\ref{eq2}) and (\ref{eq3}), and the fact that $S(N)\subseteq S(\langle E_{1}(N)\rangle)\subseteq \cdots \subseteq  S(\langle E_{n-1}(N)\rangle)$  we get equality for $i=0, 1, \cdots, n-1$. Since $S(N)=\langle E_n(N)\rangle$, $\langle E_n(N)\rangle$ is a semiprime submodule of $M$. It follows that $S(\langle E_n(N)\rangle)=\langle E_n(N)\rangle$ which completes the proof.

  \end{prf}

	\paragraph\noindent
  Let $N$ be a submodule of $M$.  If the ascending chain (\ref{chain2})   terminates at $\langle E_{n}(N) \rangle,$ then Theorem \ref{chain} says that the semiprime radical is an invariant on  the submodules $\langle E_i(N)\rangle$ for all  $i=0, 1,2, \cdots, n$. Thus, generalising Corollary \ref{cc}.

\begin{exam}\rm\label{ex}
	Let $R=\Z[X],$ $M=R\oplus R$ and $N=\{(r,s)\in M ~|~ 4r-sX \in RX^{2}\}$. $\langle E_{1}(N) \rangle=R(0,4)+XM \neq R(0,2)+XM=\langle E_{2}(N) \rangle$.
	$M$ is an $R$-module and $N$ is a submodule of $M.$ By \cite[Example 1]{Azizi} and \cite[Page 110]{PSmith}, $\langle E_{1}(N) \rangle \neq \langle E_{2}(N) \rangle$. Since $\langle E_{1}(N) \rangle$ is a preradical and the inclusion $\langle E_{1}(N) \rangle \subseteq \langle E_{2}(N) \rangle$ is strict, it follows by \cite[$\S$ VI, Proposition 1.5]{Sten} and Theorem \ref{chain} that $\langle E_{1}(N) \rangle$ is not a radical in general. The smallest radical containing the preradical $\langle E_{1}(N) \rangle$ is $S(N)$  the semiprime radical of the submodule $N$.
\end{exam}

\begin{exam}\rm The envelope functor is not right exact. Consider the canonical epimorphism $\Z\rightarrow \Z/4\Z$ of $\Z$-modules. Applying the envelope functor yields the map $0\rightarrow 2\Z/4\Z$ which is not surjective. 

\end{exam}

\paragraph\noindent 
The prime radical $\beta(R)$ of any ring $R$ coincides with $S(R)$ its semiprime radical, see \cite[Theorem 4.20]{McCoy}. This is not true in general for modules. In Corollary \ref{eq}, we give a condition under which it holds. 
\begin{cor}\label{eq}Let $N$ be a submodule of an $R$-module $M$ and $n\in \Z^{+}.$
 For modules that satisfies the radical formula of degree $n,$ $S(N)= \beta(N)$.

\end{cor}
\begin{prf}
	 The proof follows from Theorem \ref{chain} and the fact that for modules that satisfy the radical formula of degree $n,$ $\langle E_{n}(N)  \rangle=\beta(N).$
	 \end{prf}
	
\paragraph\noindent
 It then follows that the degree $n$ in Corollary \ref{eq} is a measure of how far the envelope functor is away from becoming a radical.

\paragraph\noindent Let $a\in R.$ The functor $a\Gamma_{a}:R\text{-Mod} \to R\text{-Mod}$; $M\mapsto a\Gamma_{a}(M),$ known as the locally nilradical, \cite{Ann} associates to every $R$-module $M$, a submodule $a\Gamma_{a}(M)$, where
$a\Gamma_{a}(M):=\left\{am~|~a^{k}m=0, ~m\in M, \text{for some}~k\in \Z^{+} \right\}$. It generalises the $a$-torsion functor $\Gamma_{a}$  
where $\Gamma_{a}(M):=\left\{m~|~a^{k}m=0, ~m\in M, \text{for some}~k\in \Z^{+} \right\}$. An $R$-module $M$ is {\it $a$-reduced} if $a\Gamma_{a}(M)=0$ and it is {\it reduced} if $a\Gamma_{a}(M)=0$ for all $a\in R$. For applications of both reduced modules and their generalisations, see \cite{Kim}, \cite{Ssev} and \cite{AppII}.

\begin{prop}\label{nonrad}
  For any $R$-module $M,$  $$\sum \limits_{a\in R }a\Gamma_{a}(M)= \langle E_{M}(0)  \rangle=\sum\limits_{a\in R}\underset{k}\varinjlim {\text{Hom}}_{R}(R/(a)^{k},
  aM).$$
\end{prop}
\begin{prf}
  By \cite[Proposition 4.1]{Ann}, $E_{M}(0)=\bigcup\limits_{a\in R}a\Gamma_{a}(M).$ Moreover, the smallest submodule containing $E_{M}(0)$ (resp. $\bigcup\limits_{a\in R}a\Gamma_{a}(M))$ is $\langle E_{M}(0)  \rangle$ (resp. $\sum \limits_{a\in R }a\Gamma_{a}(M)).$ Since $\Gamma_{a}$ is $R$-linear, $a\Gamma_{a}(M)=\Gamma_{a}(aM)=\underset{k}\varinjlim {\text{Hom}}_{R}(R/(a)^{k},
  aM).$  
\end{prf}

\begin{cor}
	For any ring $R$, we have 
	\begin{equation*}
		\sum \limits_{a\in R }a\Gamma_{a}(R)=\mathcal{N}(R)=\sum\limits_{a\in R}\underset{k}\varinjlim {\text{Hom}}_{R}(R/(a)^{k}, aR).
\end{equation*}\end{cor}
\begin{prf}
	Take $M=R$ in Proposition \ref{nonrad}. 
\end{prf}

\paragraph\noindent The sum of radicals need  not be a radical. Whereas $a\Gamma_{a}(M)$ is a radical, \cite[Proposition 3.1]{Ann}, 
$\langle E_{M}(0)  \rangle=\sum \limits_{a\in R }a\Gamma_{a}(M)$ is not a radical, see Example \ref{ex}.

\section{Nil modules}

\begin{defn}\rm\cite[Definition 2.2]{Beh}\label{nilpotent}
 An element $m$ of an $R$-module $ M$ is called {\it nilpotent}   if $m=\sum_{i=1}^ra_im_i$ for some $a_i\in R$, $m_i\in M$ and $r\in \Z^+$ such that $a_i^km_i=0$ ($1\leq i \leq r)$ for some $k\in \Z^+$. 
 \end{defn}
 
 \paragraph\noindent
 This definition of a nilpotent element of a module is different from the one given in \cite{Groen}.  It is superior to the one in \cite{Groen} in the sense that, unlike in \cite{Groen}, here the sum of nilpotent elements of a module is nilpotent and the collection of all nilpotent elements of a module forms a submodule.
 
 \begin{defn}\rm 
  A (sub)module is {\it nil} if every element in it is nilpotent.
  \end{defn}

 \begin{prop}\label{e}
 
 Let $M$ be an $R$-module and $a, b\in R$. 
 
 \begin{enumerate}
 
 \item Every submodule of $M$ of the form $a\Gamma_{a}(M)$ is  nil. 
 
 \item If $N_1= a\Gamma_{a}(M)$  and $N_2= b\Gamma_{b}(M)$, then $N_1+N_2$ is nil. 
 
 \item $\langle E_{M}(0)  \rangle = \sum\left\{N~:~ N~\text{is a nil submodule of}~M\right\}$ the largest nil submodule of $M$.
 
 \item A module $M$ is nil if and only if $\langle E_{M}(0)  \rangle= M$.  
 
 \item Every nil $R$-module $M$ is a sum of $R$-modules of the form $a\Gamma_a(M)$. 
 \end{enumerate}
 
   \end{prop}
   
   \begin{prf} 1) and 2) are trivial and also consequences of 3). To prove 3), we show that $\langle E_{M}(0)  \rangle$ is nil and any other submodule of $M$ which is nil is contained in it. Let $m\in \langle E_{M}(0)  \rangle$, $m=\sum_{i=1}^ra_ib_im_i$ such that $b_i^{k_i}m_i=0$ for $a_i, b_i\in R$, $m_i\in M$, $k_i\in \Z^+$ and $1\leq i \leq r$. If $k=\text{max}\{x_i\}_{i}$, then $(a_ib_i)^km_i=0$ and $m$ is nilpotent. Let $K$ be a nil submodule of $M$. For all $m\in K$, $m=\sum_{i=1}^ra_im_i$ and $a_i^km_i=0$. Then for all $1\leq i \leq r$, $a_im_i\in E_M(0)$ and $m\in \langle E_{M}(0)  \rangle$. We now prove 4). $\langle E_{M}(0)  \rangle\subseteq M$. If $M$ is nil, then by 3) $M\subseteq \langle E_{M}(0)\rangle$ and  $\langle E_{M}(0)  \rangle= M$. The converse is immediate from 3). 5) is a consequence of 4) and Proposition \ref{nonrad}.
   
   \end{prf}
   
   \paragraph\noindent
  
 A {\it torsion theory}   for  a module category $R$-Mod is a pair $(\mathcal{T}, \mathcal{F})$ of classes of $R$-modules such that: 1)
   $\text{Hom}(T, F)=0$ for all $T\in \mathcal{T}$ and $F\in \mathcal{F}$; 2) if $\text{Hom}(A, F)=0$ for all  $F\in \mathcal{F}$, then $A\in \mathcal{T}$; and 3) if $\text{Hom}(T, B)=0$ for all $T\in \mathcal{T}$,  then $B \in \mathcal{F}$.    
  $\mathcal{T}$ (resp. $\mathcal{F}$ ) is called the {\it torsion class} (resp.   {\it torsionfree class}) of the torsion theory. 
  The set of all reduced $R$-modules is  not in general a torsion-free class.   However, for rings that satisfy the radical formula, we have Proposition \ref{so}.

\begin{thm}\label{so}
Let $R$ be a ring which satisfies the radical formula. Then 
\begin{enumerate}
	\item the functor $F:R$-Mod $\to R$-Mod, $M\mapsto \langle E_{M}(0)  \rangle$ is an idempotent radical,
	
	\item the class of all reduced $R$-modules forms a torsion-free class,
	
	\item the class of all nil $R$-modules forms a torsion class.

\end{enumerate}
 
\end{thm}
\begin{prf}
  If a ring $R$   satisfies the radical formula, then $\langle E_{M}(0)  \rangle=\beta(M)$ which is a radical, i.e., the associated chain of envelopes terminates right away, i.e., $E_{2}(0)= E_{1}(0)$ if and only if $E_{1}(\langle E_{1}(0)\rangle)=E_{1}(0)$ which gives the idempotency.
	  Since $F$ is an idempotent radical, there is an associated torsion theory with a torsion-free class given by $\left \{M\in R\text{-Mod} ~ | ~ \langle  E_{{M}}({0}) \rangle= 0 \right \}$ which is the set of all reduced $R$-modules and a torsion class $\left\{M\in R\text{-Mod}~|\langle  E_{{M}}({0}) \rangle= M \right \}$ which is the set of all nil $R$-modules.
\end{prf}
 
  \paragraph\noindent
  
  Arithmetical rings \cite[Theorem 2.4]{Parkash}, Artinian rings \cite[Corollary 2.9]{Sharif}, \cite[Theorem 3.5]{Le}  and ZPI rings \cite[Corollary 2.10]{Sharif} satisfy the radical formula. A ring is called a {\it ZPI-ring} if every ideal in it is a finite product of prime ideals.

	\begin{thm}\label{rad}
		A Noetherian uniserial $R$-module $M$ satisfies the semiprime radical formula and there exists  $a\in R$ such that
		$$E_{M}(0)=a\Gamma_{a}(M)=S(M).$$  
		
	\end{thm}

	\begin{prf}	By \cite[Proposition 4.1]{Ann}, $E_{M}(0)=\bigcup\limits_{r\in R}r\Gamma_{r}(M)$.
		If $M$ is uniserial, then for any $r_{1},r_{2} \in R,$ either $r_{1}\Gamma_{r_{1}}(M)\subseteq r_{2}\Gamma_{r_{2}}(M)$ or $r_{2}\Gamma_{r_{2}}(M)\subseteq r_{1}\Gamma_{r_{1}}(M).$  Without loss of generality, suppose that $r_{1}\Gamma_{r_{1}}(M)\subseteq r_{2}\Gamma_{r_{2}}(M)$. So, we can form a chain $r_{1}\Gamma_{r_{1}}(M)\subseteq r_{2}\Gamma_{r_{2}}(M) \subseteq r_{3}\Gamma_{r_{3}}(M) \subseteq \cdots$ which stabilizes since $M$ is Noetherian. So, $\bigcup\limits_{r_{i}\in R}r_{i}\Gamma_{r_{i}}(M)=r_{k}\Gamma_{r_{k}}(M)$ for some $r_{k}\in R$  and $\langle E_M(0)\rangle=E_{M}(0)=r_{k}\Gamma_{r_{k}}(M).$ Taking $a=r_{k}$ gives the desired result. Since $a\Gamma_{a}(M)$ is a radical and $\langle E_M(0)\rangle$ is in general a preradical,  the ascending chain $0 \subseteq \langle E_M(0)\rangle \subseteq \langle E_2(0)\rangle\subseteq \langle E_3(0)\rangle \subseteq \cdots$ stabilizes at a radical \cite[Proposition 1.5, Chapter VI]{Sten}. Furthermore, by Corollary \ref{cc}, this radical is the semiprime radical $ S(M)$ of $M.$ Hence, $E_{M}(0)=a\Gamma_{a}(M)=S(M).$  Since $M$ is uniserial and Noetherian, for any submodule $N$ of $M$, $M/N$ is also uniserial and Noetherian. As before (for the case when $N=0$), we have for some $a\in R$,
		
		$$E_{\frac{M}{N}}(\bar{0})=\frac{E_M(N)}{N}=\bigcup_{r_i\in R}r_i\Gamma_{r_i}\Bigl(\frac{M}{N}\Bigr)=a\Gamma_a\Bigl(\frac{M}{N}\Bigr)=S\Bigl(\frac{M}{N}\Bigr)=\frac{S(N)}{N}.$$ It follows that $E_M(N)=S(N)$ for any submodule $N$ of $M$.
			\end{prf}

	\begin{cor}\label{Not}
		If $R$ is a Noetherian uniserial ring, then there exists $a\in R$ such that $$\mathcal{N}(R)=a\Gamma_{a}(R).$$ 
	\end{cor}
	
	\paragraph\noindent
	By Proposition \ref{e}(3),  the submodules $\beta(M)$ and $S(M)$ of $M$ are not nil in general since  $\langle E_{M}(0)  \rangle$ is the largest nil submodule $M$ and it is in general strictly contained in $S(M)\subseteq \beta(M)$. However, we have Proposition \ref{nil}.
	
	\begin{prop} \label{nil}	
	The following statements hold: 
	\begin{enumerate}
	
	\item If $R$ satisfies the radical formula, then for every  $R$-module $M$, the prime radical $\beta(M)$ of $M$ is a nil submodule.
	
	\item If $M$ is a Noetherian uniserial  $R$-module, then $S(M)$ is a nil submodule of $M$.	
	\end{enumerate}
	
	\end{prop}
		
 \begin{prf} If $R$ satisfies the radical formula, then $\beta(M)= \langle E_M(0)\rangle$ and by Proposition \ref{e}, $ \langle E_M(0)\rangle$ is nil. Part 2 follows from Theorem \ref{rad} and the fact that  $a\Gamma_a(M)$ is nil.
 
 \end{prf}

  \begin{prop} For any $R$-module $M$ and $a\in R$,
  \begin{enumerate}
  \item $a\Gamma_a(M)\cong \Gamma_a(M)/(0:_{\Gamma_a(M)}a)$,
  
  \item $M$  is $a$-reduced if and only if $(0:_{\Gamma_a(M)}a)=\Gamma_a(M)$.
  
  \item  $(0:_{\Gamma_a(M)}a)=0$  if and only if $\Gamma_a(M)\cong a\Gamma_a(M)$.
  
  \end{enumerate}

  \end{prop}

  \begin{prf} Define a map $\gamma: \Gamma_a(M)\rightarrow \Gamma_a(M)$ by $\gamma(m)=am$. $\text{im}\gamma = a\Gamma_a(M)$ and $\text{ker}\gamma =  (0:_{\Gamma_a(M)}a)$. 1) follows from the isomorphism theorem. From 1), it is evident that $a\Gamma_a(M)=0$ if and only if $(0:_{\Gamma_a(M)}a)=\Gamma_a(M)$ which leads to 2); and $a\Gamma_a(M)\cong\Gamma_a(M)$ if and only if $(0:_{\Gamma_a(M)}a)=0$ which is 3).

  \end{prf}
  
  \begin{rem}\rm 
  For any $R$-module $N$, it is  well known that there exists an injective $R$-module $E(N)$ called the injective envelope or injective hull of $N$. It is natural for one to ask whether there is a relationship between  this notion  and the envelope studied in this paper. In general, these two notions are different. Whereas, if $N$ is a submodule of $M$, $E(N)$  just like $E_M(N)$ both contain $N$, (i.e., both envelope $N$), $E(N)$ is always an injective $R$-module but $E_M(N)$ need not even be a module in general. However, for $R=N=\Z$, the ring of integers, we have $$E(\Z)=\Q=E_{\Q}(\Z).$$ The first equality is well known. For the second  equality, first note that in general, $E_{\Q}(\Z)\subseteq \Q$. For the reverse inclusion, take $m\in\Q$, $m=\frac{a}{b}$, where $a, b\in \Z$. It follows that $m=b(\frac{a}{b^2})$,  $b^2(\frac{a}{b^2})=a\in \Z$  and $m\in E_{\Q}(\Z)$ such that $E_{\Q}(\Z)=\Q$ as required.
   
  \end{rem}

{\bf{Acknowledgement}}

\paragraph\noindent
We acknowledge support from the Makerere University Research and Innovation Fund (RIF) and from the International Science Programme (ISP) through the Eastern Africa Algebra Research Group (EAALG). We are grateful to the referee for the comments that  improved the presentation of the paper.\\

{\bf{Disclosure Statement}}: The authors declare that they have no conflict of interest.

\end{document}